\documentclass{proc}
\usepackage[utf8]{inputenc}
\usepackage[centertags]{ amsmath}

\usepackage{amssymb,amsthm}
\usepackage{amsfonts}
\usepackage[colorlinks=true,linkcolor=blue,citecolor=blue]{hyperref}
\usepackage{url}
\usepackage{multirow}
\usepackage{csquotes}
\usepackage{subcaption}
\usepackage[rgb]{xcolor}
\usepackage{tikz}
\usetikzlibrary{decorations.pathmorphing}
\usepackage{adjustbox}

\usepackage{graphicx}

\title{Phase Portraits of Hyperbolic Geometry}
\author{Scott B. Lindstrom \\ CARMA\\ University of Newcastle \and Paul Vrbik \\  University of Toronto Mississauga}

\newtheorem{theorem}{Theorem}

\theoremstyle{definition}

\begin{document}

\maketitle

\begin{displayquote}
	Sometimes it is easier to see than to say. --- Jonathan M. Borwein
\end{displayquote}
Borwein was quite fond of saying this, and it seems fitting that his first posthumously published book was \emph{Tools and Mathematics: Instruments for Learning} with Luc Trouche and John Monaghan. In his chapter \cite{Tools}, he included, along with his own commentary, a quote he particularly liked:
\begin{displayquote}
	Long before current graphic, visualisation and geometric tools were available,
	John E. Littlewood, 1885-1977, wrote in his delightful Miscellany:
	\begin{displayquote}
	A heavy warning used to be given [by lecturers] that pictures are not rigorous; this has never had its bluff called and has permanently frightened its victims into playing for safety. Some pictures, of course, are not rigorous, but I should say most are (and I use them whenever possible myself).\cite[p.53]{Littlewood}
	\end{displayquote}
\end{displayquote}
The popularity of Arnold and Rogness' video \emph{M\"{o}bius transformations revealed} \cite{ARvideo} serves to highlight that not only may pictures be rigorous, they may appeal to a much broader audience than researchers, a point they take care to emphasize in their follow-up article \cite{AR}.

Tristan Needham's award-winning \emph{Visual Complex Analysis} \cite{Needham} and Elias Wegert's \cite{Wegert} \emph{Complex Functions: An Introduction with Phase Portraits} have employed visual tools to reorient the learning of complex analysis around the geometric intuition which, historically, motivated important contributions to the field. In particular, Wegert's book and his article with Semmler \cite{WS} have done much to draw attention to \emph{phase portraits} which plot complex functions $f:\mathbb{C}\rightarrow \mathbb{C}$ by assigning colors to points according to the \emph{phase} of their image under $f$. One such assignment is the map to the color wheel given by
\begin{align*}
\mathcal{C}_{\mathbb{C}}&:\mathbb{C}\rightarrow [0,2\pi] \quad \text{by}\quad z \mapsto \arg \circ (f(z)/|f(z)|).
\end{align*}

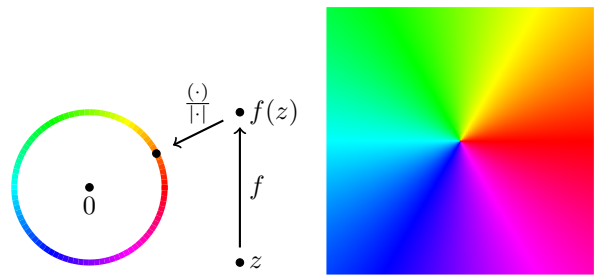
\begin{figure}[h]
	\begin{center}
		\hspace*{\fill}
		\begin{tikzpicture}

		%draw and label z
		\draw [fill] (2,-1) circle [radius=0.05];
		\node [right] at (2,-1) {$z$};
		
		%draw and label f(z)
		\draw [fill] (2,1) circle [radius=0.05];
		\node [right] at (2,1) {$f(z)$};
		
		%draw arrow from z to f(z) and label
		\draw [thick,->] (2,-.8) -- (2,.8); 
		\node [right] at (2,0) {$f$};
		
		%draw final arrow to point on unit circle
		\draw [thick,->] (1.78,0.89) -- (1.11,0.56);
		\node [above] at (1.44,0.72) {$\frac{(\cdot)}{|\cdot|}$};

		%draw and albel the origin
		\draw [fill] (0,0) circle [radius=0.05];
		\node [below] at (0,0) {$0$};
		
		%Here I plot my unit circle
		%What I'm actually doing is plotting a bunch of 
		%quadrilaterals, each with its own color. I plot enough of them and I get the circle.
		\foreach \x in {0,0.0111,...,1} {
			\definecolor{currentcolor}{hsb}{\x, 1, 1}
			\draw[draw=none, fill=currentcolor]
			(360*\x - 2 : 0.96) -- (360*\x - 2 :1.04)
			-- (360*\x + 2:1.04) -- (360*\x + 2 :0.96) -- cycle;
		}
		
		%draw final point on unit circle
		\draw [fill] (.89,0.45) circle [radius=0.05];
		\end{tikzpicture}
		\hspace*{\fill}\includegraphics[width=.45\linewidth]{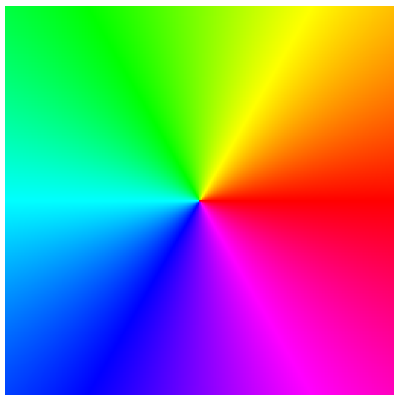}\hspace*{\fill}
	\end{center}
	\caption{Construction of a complex phase portrait}\label{fig:complexcase}	
\end{figure}

\begin{figure}[h]
	\hspace*{\fill}\includegraphics[width=.8\linewidth]{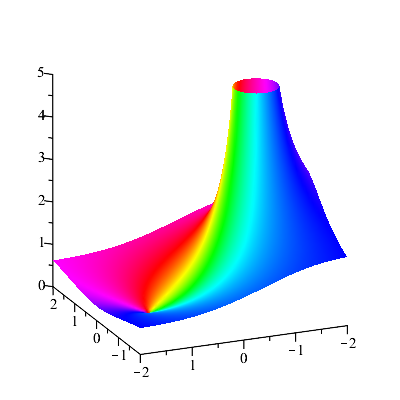}\hspace*{\fill}\caption{A 3d phase portrait with modulus.}\label{fig:complexmodular}
\end{figure}

The construction is shown at left in Figure~\ref{fig:complexcase} while a plot of the identity is at right. The modulus $|f(z)|$ may be represented by contour lines or by the height of a surface on a 3D plot as in Figure~\ref{fig:complexmodular}.

We introduce a new method for visualizing hyperbolic geometry by defining color maps on surfaces commonly employed to understand hyperbolic space. To do so we exploit the connections between hyperbolic geometry and the complex plane explored by Needham in \emph{Visual Complex Analysis} \cite{Needham}. 

The reinvention is nontrivial: the classical method of simply applying a phase portrait following a conformal (angle-preserving) map does not extend in a useful way. This forces us to revisit the original motivations of phase plotting in order to find a suitable method of applying color maps in the context of a new geometry.

Wegert and Semmler note that credit for the invention of phase plots cannot be given to a single person \cite{WS}. Phase portraits were predated by \emph{analytic landscapes} (3d modulus plots originally in grayscale). Combining the modulus plot (which gives unique heights to points lying on concentric circles about $0$) with the phase plotting (which assigns unique colors to the preimages of half lines intersecting these circles) the viewer effectively obtains a four dimensional plot.

Taking our cues from history, we adopt the principle of coloring hyperbolic ``lines''---in our setting geodesics---a single color. In so doing, we uncover natural connections with the classical idea of phase plotting and discover a beautiful way to visualize hyperbolic geometry. 

This work is inspired by Jonathan Borwein's chapter, ``The Life of Modern Homo Habilis Mathematics: Experimental Computation and Visual Theorems,'' \cite{Tools}, by his interest in the phase portraits which included a joint contribution with Armin Straub to Wegert's 2016 \emph{Complex Beauties Calendar}\cite{calendar}, and by the aforementioned seminal works in complex analysis \cite{Needham}\cite{Wegert}.

\section{Non-Euclidean Geometry}

\begin{figure}[h]
	\begin{center}
		\includegraphics[width=.62\linewidth]{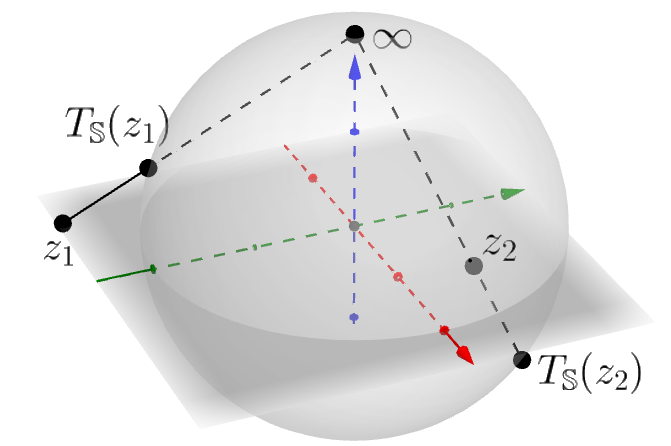}\hspace*{\fill}
		\adjustbox{trim={.21\width} {.21\height} {.21\width} {.21\height},clip}%
		{\includegraphics[width=.69\linewidth]{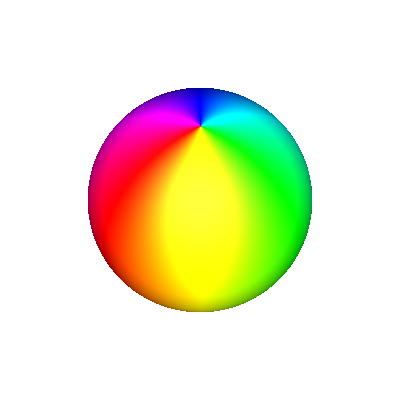}}
	\end{center}
	\caption{Using stereographic projection to plot a phase portrait for a direct motion on the Riemann sphere.}\label{fig:Riemann}
\end{figure}

Throughout, ${\rm Re}$ is the real part map $x+yi\mapsto x$ on $\mathbb{C}$, ${\rm Im}$ is the imaginary part map $x+yi \mapsto y$ on $\mathbb{C}$, and the plane is taken to be $\mathbb{C}$.

\emph{Conformal} mappings preserve the angles at which curves meet, and \emph{anticonformal} mappings preserve the magnitude of those angles while reversing their orientation. Motions preserve distances between points; direct motions are conformal while opposite motions are anticonformal. We use the term \emph{motion} instead of \emph{isometry} in order to be consistent with Needham, to whom we refer the reader for any additional details about these definitions or the preliminaries which follow \cite[pp.~31--36, Ch.~6]{Needham}.

Any direct motion is just the composition $\mathcal{R}_{L_2}\circ \mathcal{R}_{L_1}$ of two reflections $\mathcal{R}_{L_2}, \mathcal{R}_{L_1}$ in the lines $L_2,L_1$. In $\mathbb{C}$, if $L_1,L_2$ are parallel, the direct motion is a translation. If they intersect, it is a rotation. Reflecting across a line in $\mathbb{C}$ is an opposite motion whereas scaling is conformal but not a direct motion. In $\mathbb{C}$, any direct motion is affine and of the form $z\mapsto e^{i\theta}z+b$. 

The parallel axiom states that for any line $L$ and point $p\notin L$, there exists a unique line $L'$ parallel to $L$ such that $p\in L'$. There are two alternatives to the parallel axiom:
\begin{enumerate}
	\item \emph{Spherical axiom:} there is no line through $p$ parallel to $L$.
	\item \emph{Hyperbolic axiom:} there are multiple lines through $p$ parallel to $L$.
\end{enumerate}

We may understand non-Euclidean geometries by taking a ``line'' to simply be a \emph{geodesic}. A \emph{geodesic segment} is a curve segment of shortest length connecting two points; we take a \emph{geodesic} to be a curve whose segments are geodesic segments and which is contained in no other curve with this property. On a sphere any such ``line'' divides the sphere into two hemispheres (i.e. a \emph{great circle}) and any two ``spherical lines'' will intersect. Thus the geometry of a sphere obeys the spherical axiom.

The sphere also has \emph{constant curvature}: any piece of the sphere cut out and grafted onto the sphere elsewhere fits the surface without needing to be stretched. We say of such surfaces that their \emph{intrinsic} geometry is everywhere the same. Consider a carpet unrolled onto $\mathbb{C}$ without being stretched in any way. Though the extrinsic geometry of the carpet has changed, its intrinsic geometry has not because it is already the same as that of $\mathbb{C}$.

The Riemann sphere $\mathbb{S}$ and $\mathbb{C}$ are connected via the conformal map of stereographic projection $T_\mathbb{S}:\mathbb{C}\rightarrow \mathbb{S}$ and its inverse $T_\mathbb{S}^{-1}:S\rightarrow \mathbb{C}$. Stereographic projection establishes a correspondence between each point in the plane and a point on the sphere with its north pole removed, as may be seen in Figure~\ref{fig:Riemann}, where we have labeled the north pole $\infty$ for convenience. Geometrically, we construct the map by taking a sphere of radius $1$ and centering it at the origin $0$ in the complex plane; to assign a point $z$ in the plane to a point $\widehat{z}$ on the sphere, we then draw the line connecting the north pole of the sphere to $z$. The point $\widehat{z}$ is the unique point where the line intersects the sphere. Where $z := x+iy$ and $\widehat{z}:= (X,Y,Z)$, the map is given explicitly by
\begin{eqnarray*}
	T_\mathbb{S}^{-1}: &  \quad x+iy &= \frac{X+iY}{1-Z},\\
	T_\mathbb{S}:  &\quad X+iY &= \frac{2z}{1+|z|^2},\\
	& \quad \text{and} \quad Z &= \frac{|z|^2-1}{|z|^2+1},
\end{eqnarray*}
the derivation of which may be found in \cite[p. 146]{Needham}. With this construction, it is clear why we label the north pole $\infty$. As the modulus of a point $z$ approaches infinity, its counterpart $\widehat{z}$ approaches the north pole of the sphere. 

Stereographic projection formed the basis of the visualizations employed in the \emph{M\"{o}bius Transformations Revealed} \cite{ARvideo} and the corresponding article \cite{AR}. Phase plotting on the sphere itself has already been explored in depth by Wegert and Semmler in their article and later by Wegert in his book \cite{WS}\cite{Wegert}. 

A direct motion on $\mathbb{S}$, just like any direct motion in $\mathbb{C}$, may be expressed as the composition of two ``reflections'' about ``lines'' $C_2,C_1$ on the sphere and computed in $\mathbb{C}$ by
\begin{equation*}
T_\mathbb{S} \circ \mathcal{R}_{L_2} \circ \mathcal{R}_{L_1} \circ T_\mathbb{S}^{-1}
\end{equation*}
where $L_2:=T_\mathbb{S}^{-1}(C_2)$ and $L_1:=T_\mathbb{S}^{-1}(C_1)$ are regarded as the representatives in $\mathbb{C}$ of the sphere ``lines'' $C_2,C_1$. We adopt this approach of referring to lines $C_j$ in non-Euclidean space by their representatives $L_j$ in the plane. Figure~\ref{fig:Riemann} illustrates the conformal map $T_\mathbb{S}$ at left and a direct motion (the identity) at right. On the sphere all lines intersect and so the sphere possesses only one direct motion: rotation.

\section{Hyperbolic Space}

\begin{figure}[h]
	\begin{center}
		\includegraphics[width=\linewidth]{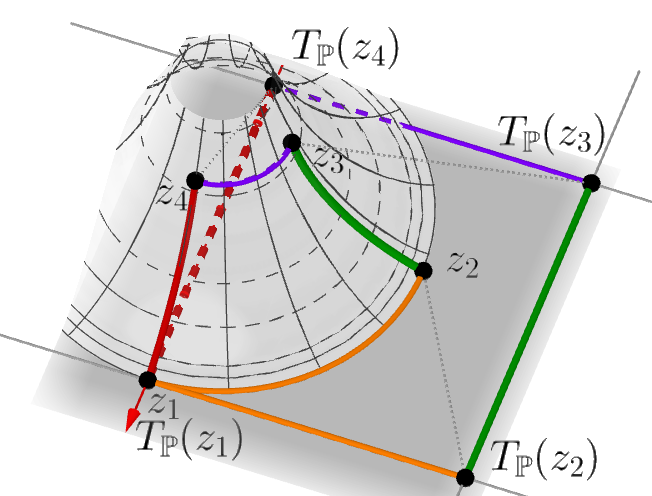}
	\end{center}
	\caption{The conformal map $T_{\mathbb{P}}$ from the pseudosphere to the upper half plane.}\label{fig:Pseudosphere1}
\end{figure}

\subsection{Pseudosphere}

The pseudosphere ($\mathbb{P}$) is a surface with constant curvature satisfying the axioms of  hyperbolic geometry. The physical pseudosphere is made by rotating a curve called a \emph{tractrix} through $2\pi$ about an axis; in terms of Figure~\ref{fig:Pseudosphere1}, the curve along the pseudosphere that connects points $z_1$ and $z_4$ is one example of a tractrix, while the curve connecting $z_2$ and $z_3$ is another. 

A conformal mapping exists from the pseudosphere to $\mathbb{C}_+$, illustrated at left in Figure~\ref{fig:Pseudosphere1}. The pseudosphere has been superimposed on the complex plane with base centered at $(0,2i)$ for the sake of readability. Each point on the surface is labeled with a coordinate  $(\theta, \sigma)$ where $\theta$ corresponds to an angle of rotation about the center axis of the pseudosphere and $\sigma$ is the distance traveled along the tractrix generator from the rim of the sphere. With this labeling system, points along the same tractrix in Figure~\ref{fig:Pseudosphere1} have the same $\theta$ value whereas points along the same level curve have the same $\sigma$ value. We have chosen to illustrate with points
$$z_1=(0,0),\;\; z_2=(\pi/2,0),\;\; z_3=(\pi/2,1),\;\; z_4=(0,1).$$ The conformal map is
 \begin{equation*}
 \mathbb{P}\rightarrow [0,2\pi)\times[1,\infty) \quad \text{by} \quad (\theta,\sigma) \mapsto (\theta,e^\sigma i) =: (x,y),
 \end{equation*}
where the line segment $[0,2\pi)\times \{1\}$ is the image of the pseudosphere's rim. For two infinitesimally close neighbors $x+yi$ and $(x+dx)+(y+dy)i$, the metric associated with this map is
\begin{equation*}
d\widehat{s} = \frac{ds}{y} = \frac{\sqrt{dx^2+dy^2}}{y}.
\end{equation*}
See \cite[p. 298]{Needham}.

The representation of hyperbolic space with a pseudosphere carries two main disadvantages when compared with the hyperbolic plane model $\mathbb{C}_+$.
\begin{enumerate}
	\item The first is that, in the plane, closed loops can be continuously shrunk to a point.
	\item The second is that, in the plane, geodesics may be extended infinitely far in either direction.
\end{enumerate}
For the former problem, Beltrami proposed the remedy of imagining the pseudosphere to be a roll of paper wrapping around infinitely both clockwise and counterclockwise. Under this construction, when we look at the pseudosphere, we are only \emph{seeing} one rotation of the wrapped paper rather than the entire surface. We may see more if we take the visual liberty of ``unrolling'' the paper as in the construction of Dini's surface, though this carries the cost of sacrificing the extrinsic geometry. 

Beltrami's wrapping interpretation allows us to extend the map as follows: 
\begin{equation*}
\mathbb{P}\rightarrow (-\infty,\infty)\times[1,\infty) \quad \text{by} \quad (\theta,\sigma) \mapsto (\theta,e^\sigma i).
\end{equation*}

For the second disadvantage, our best remedy is simply to imagine that our surface extends ``down'' beyond the rim of the pseudosphere into ``negative'' tractrix distances. This allows us we recover all of hyperbolic space $\mathbb{C}_+$. The complete map is then given by
\begin{equation*}
T_{\mathbb{P}}:\mathbb{P}\rightarrow \mathbb{C}_+ \quad \text{by} \quad (\theta,\sigma) \mapsto (\theta,e^\sigma i).
\end{equation*}

The map sends tractrices to vertical lines in $\mathbb{C}_+$ and all other hyperbolic lines to semi-circles centered on the real axis; this fact will prove very important for the construction of our coloring assignments. For additional history and a more detailed construction of the map and the above exposition, see \cite[pp. 296--299]{Needham}.

\begin{figure}[ht]
	\begin{center}
		\hspace*{\fill}\includegraphics[width=.46\linewidth]{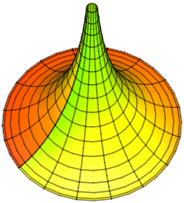}\hspace*{\fill}\includegraphics[width=.46\linewidth]{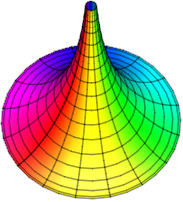} \hspace*{\fill}
	\end{center}
	\caption{Left: a naive attempt at phase plotting on a pseudosphere. Right: uniquely coloring the tractrix lines of a pseudosphere.}\label{fig:Pseudosphere}
\end{figure}
The natural temptation is to plot functions on the pseudosphere by mapping to $\mathbb{C}$, applying the desired function, and coloring according to the argument of the output. This color map
\begin{equation*}
\arg \circ f \circ T_{\mathbb{P}}
\end{equation*}
is the natural analog of the color map for phase plotting on the sphere
\begin{equation*}
\arg \circ f \circ T_{\mathbb{S}}.
\end{equation*}
However, this is ultimately useless as can be seen at left in Figure~\ref{fig:Pseudosphere}. At a glance it becomes clear that the space constraints of the pseudosphere preclude much useful information from being shown. However, the true fatal flaw of this method is that unique colors no longer show the preimages of geodesics.

Clearly a better way is needed, and   fortunately one presents itself. If we set aside the complex plane for a moment and think instead about how to apply the color spectrum to the pseudosphere, one method stands out: we may associate the rim of the pseudosphere with the color wheel. If we go one step further and color each point on the pseudosphere the same color as its nearest point along the rim---as at right in Figure~\ref{fig:Pseudosphere}---we obtain a method which \emph{does} assign unique colors to unique geodesics! Particularly, these geodesics are the tractrix generators with $\theta$ coordinates in $[0,2\pi)$.

The method still has its limitations: It does not provide colors for points which are mapped ``out of view'' to points with $\theta$ coordinates outside of $[0,2\pi)$. We color those points black. The complete map is
\begin{align*}
\mathcal{C}_{\mathbb{P}}:\mathbb{P}&\rightarrow [0,2\pi]\cup\{\text{black}\} \\
\text{by} \quad z &\mapsto \begin{cases} {\rm Re}  \circ f\circ T_{\mathbb{P}}(z) & {\rm Re}  \circ f\circ T_{\mathbb{P}}(z) \in [0,2\pi] \\
\text{black} & \text{otherwise} \end{cases}.
\end{align*}

Any direct motion in hyperbolic space is merely the composition of two reflections $\mathcal{R}_{C_2} \circ \mathcal{R}_{C_1}$ in h-lines $C_1,C_2$ which we identify by their representatives $L_1:=T_\mathbb{D}^{-1}(C_1)$ and  $L_2:=T_\mathbb{D}^{-1}(C_2)$ in $\mathbb{C}_+$. Note that ``reflection'' about a line $C_j$ corresponding to a circle $L_j$ in $\mathbb{C}_+$ is just inversion in $L_j$. For hyperbolic space, there are \emph{three kinds} of direct motion: 
\begin{enumerate}
	\item \emph{Rotations} occur when $L_1,L_2$ intersect.
	\item \emph{Translations} occur when $L_1,L_2$ are ``ultra-parallel'' (they do not intersect, even at infinity).
	\item \emph{Limit rotations} occur when $L_1,L_2$ are \emph{asymptotic}, meaning they intersect at infinity (which includes the real axis in $\mathbb{C}_+$).
\end{enumerate}

\begin{figure}[h]
	\begin{center}
	\includegraphics[width=.9\linewidth]{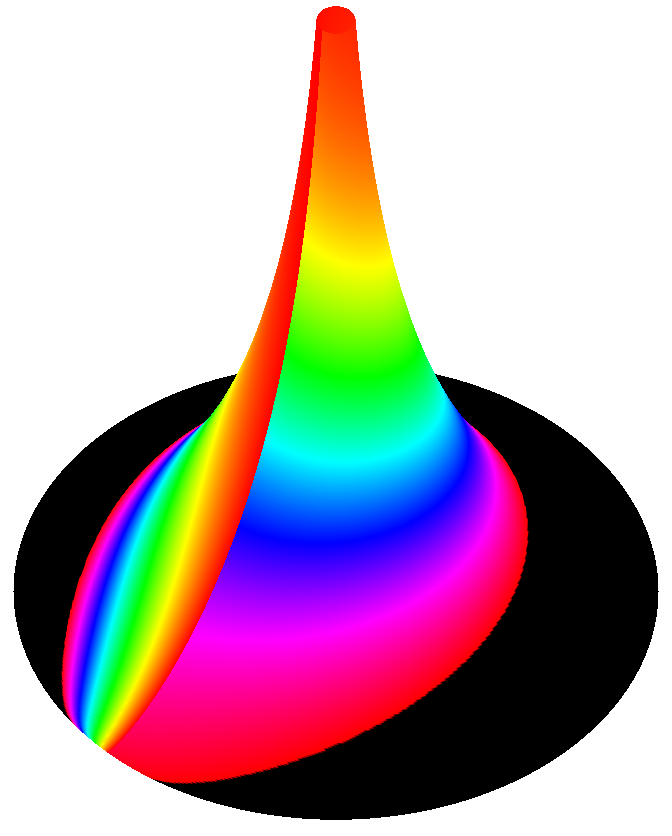}
	\end{center}
	\caption{An h-rotation plotted on a pseudosphere.} \label{fig:pseudo}
\end{figure}

Equipped with our first coloring view into hyperbolic space, we set out to visualize our first direct motion: a rotation. Let $L_1$ and $L_2$ correspond to circles of radius $2\pi$ centered at $0$ and $2\pi$ respectively. The M\"{o}bius transformation is $f(z)=4\pi^2/(2\pi-z)$ in $\mathbb{C}_+$ and the portrait is shown in Figure~\ref{fig:pseudo}. For direct motions, the preimages of geodesics are \emph{themselves} geodesics, and so we see other ``lines'' appear. 

\begin{figure}[h]
	\begin{center}
	\includegraphics[width=.9\linewidth]{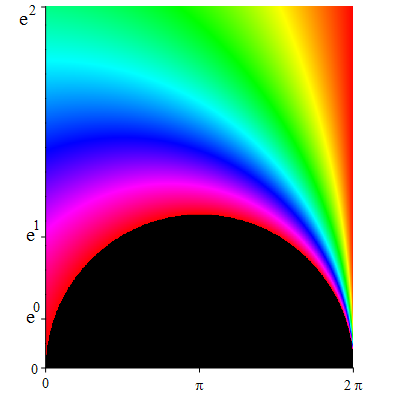}\end{center}
	\caption{An h-rotation plotted in $\mathbb{C}_+$.} \label{fig:pseudo2}
\end{figure}

Peering into $\mathbb{C}_+$ in Figure~\ref{fig:pseudo2}, we glean additional insights about what we are seeing on the pseudosphere, allowing us to explain the behavior. The preimages of tractrix lines are semicircles in $\mathbb{C}_+$ which intersect at the infinity corresponding to $2\pi$ in $\mathbb{C}_+$. The preimages of those tractrix lines we see on the pseudosphere are those ``outside'' the semicircle of radius $\pi$ centered at $\pi$. Semi-circles ''inside'' are sent to tractrix lines to the right of $2\pi$. The preimages of tractrix lines to the left of $0$ are semicircles with centers to the right of $2\pi$ meeting at $2\pi$.

Having seen what a direct motion looks like on the pseudosphere, we would really like to have a way of viewing \emph{all} of hyperbolic space. Fortunately, such representations exist.

\subsection{Poincar\'{e} Disc}

\begin{figure}[h]
	\begin{center}
		\begin{tikzpicture}[scale=2.0]
		
		%draw the semicircular arc for K
		\draw [gray,domain=0:180] plot ({1.41*cos(\x)}, {-1+1.41*sin(\x)});
		\node [above right, gray] at (1.41,-1) {$K$};
		
		%draw and label the unit circle
		\draw [thick] (0,0) circle [radius=1];
		
		%draw and label the point at i
		\draw [fill,gray] (0,-1) circle [radius=0.04];
		\node [above] at (0,-1) {$-i$};
		
		%draw and label the origin
		\draw [fill,gray] (0,0) circle [radius=0.04];
		\node [below,gray] at (0,0) {$0$};
		
		%draw and label the point 1
		\draw [fill,gray] (-1,0) circle [radius=0.04];
		\node [above left,gray] at (-1,0) {$-1$};
		
		%draw the axis
		\draw (-1.41,0) -- (2,0);
		
		%horizontal reds
		\foreach \y in {0,0.1,...,0.5} {
			\draw[red] (1,\y) -- (2,\y);
			
			%\draw [red,domain=1:2] plot ({(4*\x)/(2*|\x|^2+2*|\y+1|^2)}, {4*(\y+1)/(2*|x|^2+2*|\y+1|^2)-1});
			
			\draw [red,domain=1:2] plot ({(4*\x)/(2*\x*\x+2*(\y+1)*(\y+1))}, {4*(\y+1)/(2*\x*\x+2*(\y+1)*(\y+1))-1});
			
		}
		
		%vertical reds
		\foreach \y in {1,1.1,...,2.1} {
			\draw[red] (\y,0) -- (\y,0.4);
			
			%\draw [red,domain=1:2] plot ({(4*\y)/(2*|\y|^2+2*|\x+1|^2)}, {4*(\x+1)/(2*|x|^2+2*|\x+1|^2)-1});
			
			\draw [red,domain=0:0.4] plot ({(4*\y)/(2*\y*\y+2*(\x+1)*(\x+1))}, {4*(\x+1)/(2*\y*\y+2*(\x+1)*(\x+1))-1});
			
		}
		
		%horizontal blues
		\foreach \y in {0,0.1,...,0.4} {
			\draw[blue] (-.5,\y) -- (0,\y);
			
			%\draw [red,domain=1:2] plot ({(4*\x)/(2*|\x|^2+2*|\y+1|^2)}, {4*(\y+1)/(2*|x|^2+2*|\y+1|^2)-1});
			
			\draw [blue,domain=-.5:0] plot ({(4*\x)/(2*\x*\x+2*(\y+1)*(\y+1))}, {4*(\y+1)/(2*\x*\x+2*(\y+1)*(\y+1))-1});
			
		}
		
		%vertical blues
		\foreach \y in {-.5,-.4,...,0.1} {
			\draw[blue] (\y,0) -- (\y,0.3);
			
			%\draw [red,domain=1:2] plot ({(4*\y)/(2*|\y|^2+2*|\x+1|^2)}, {4*(\x+1)/(2*|x|^2+2*|\x+1|^2)-1});
			
			\draw [blue,domain=0:0.3] plot ({(4*\y)/(2*\y*\y+2*(\x+1)*(\x+1))}, {4*(\x+1)/(2*\y*\y+2*(\x+1)*(\x+1))-1});
			
		}
		
		\end{tikzpicture}	
	\end{center}\caption{Anticonformal map $I_{K}$}\label{fig:anticonformaldisc}
\end{figure}
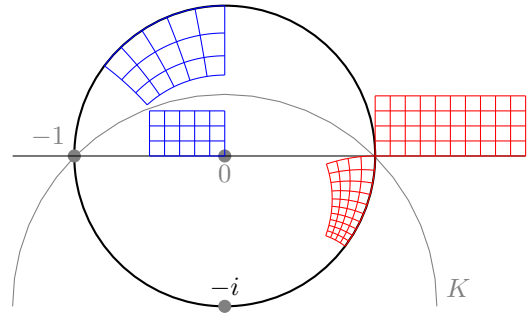

Hyperbolic space may also be visualized using what has come to be known as the Poincar\'{e} disc $\mathbb{D}$, a construction described by Beltrami. A conformal map $T_{\mathbb{D}}$ sends $\mathbb{C}_+$ to the unit disc in $\mathbb{C}$ by inverting in a circle $K$ of radius $\sqrt{2}$ centered at $-i$ and then conjugating. If we first conjugate and then invert in the circle, we obtain $T_{\mathbb{D}}^{-1}$. The map $I_K$ of inversion in the circle $K$ is shown in Figure~\ref{fig:anticonformaldisc}, and its 3d phase plot is shown in Figure~\ref{fig:complexmodular}.
%\begin{equation*}
%T_{\mathbb{D}}:\mathbb{C}_+\rightarrow \mathbb{D} \quad \text{by} \quad z \rightarrow \frac 
%{iz+1}{z+i}\end{equation*}

In $\mathbb{D}$ any hyperbolic line $C_j$ appears as a circular arc meeting the edge of $\mathbb{D}$ at a ninety degree angle or as a diameter dividing the disc in half (which may be thought of as a semicircular arc of infinite radius). As with the pseudosphere, under this map the $L_j$ representing the $C_j$ are vertical half lines and semicircles centered on the real axis in $\mathbb{C}_+$.

Our method on the pseudosphere of coloring tractrix lines is the motivation for our first attempt at phase plotting on $\mathbb{D}$. As before, we may uniquely identify any tractrix line by the $\theta$ coordinate of all its points, and the tractrix lines correspond to the vertical half-lines in $\mathbb{C}_+$ which may be uniquely identified by the real parts of all of their points. 

Thus the map ${\rm Re} \circ f \circ T_{\mathbb{D}}^{-1}$ will take the preimages of unique tractrices to unique real numbers. We need only to define a map from $\mathbb{R}$ to $(0,2\pi)$ in order to have a color wheel assignment which assigns unique colors to the preimages of unique tractrix lines. In fact, looking in Figure~\ref{fig:anticonformaldisc}, we realize that we already have a map which takes $\mathbb{R}$ to the unit circle: $I_K$. We may use $T_{\mathbb{D}}$ as our map and treat the unit circle as the color wheel!

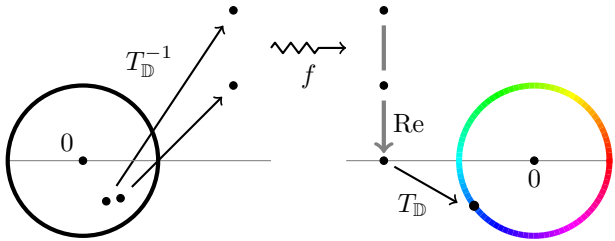
\begin{figure}[ht]
	\begin{tikzpicture}
	
	%draw the preimage points
	\draw [fill] (-5.5,-0.5) circle [radius=0.05];
	\draw [fill] (-5.69,-0.54) circle [radius=0.05];
	
	%draw the points T(z) and T(z2)
	\draw [fill] (-4,2) circle [radius=0.05];
	\draw [fill] (-4,1) circle [radius=0.05];
	
	%draw arrows connecting preimage points and image points
	\draw[->,thick] (-5.35,-0.35) -- (-4.15,0.85);
	\draw[->,thick] (-5.5548,-0.3368) -- (-4.1352,1.7968);
	\node [above left] at (-4.676,0.984) {$T_{\mathbb{D}}^{-1}$};
	
	%draw the squiggly function arrow
	\draw [->,thick,
	line join=round,
	decorate, decoration={
		zigzag,
		segment length=6,
		amplitude=2,post=lineto,
		post length=10pt				
	}] (-3.5,1.5) -- (-2.5,1.5);
	
	%label the squiggly function arrow
	\node [below] at (-3,1.4) {$f$}; 
	
	%draw the real axes
	\draw [gray] (-2.5,0) -- (1,0);
	\draw [gray] (-7,0) -- (-3.5,0);
	
	%draw points for image points of f
	\draw [fill] (-2,2) circle [radius=0.05];
	%\node [right] at (0,0) {$0$};
	\draw [fill] (-2,1) circle [radius=0.05];
	%\node [below] at (0,0) {$0$};
	
	%draws the real part arrows down
	\draw [ultra thick, gray,->] (-2,.8) -- (-2,0.1) ;
	\node [right] at (-2,.5) {Re};
	\draw [ultra thick, gray] (-2,1.8) -- (-2,1.2) ;
	
	%final arrow to the color wheel
	\draw [thick, black,->] (-1.86,-0.08) -- (-1.02,-.56) ;
	\node [below right] at (-1.94,-0.32) {$T_{\mathbb{D}}$};

	%Draw the left circle and label it
	\draw [ultra thick] (-6,0) circle [radius=1];
	
	%draw Real value point and label
	\draw [fill] (-2,0) circle [radius=0.05];
	%\node [below left] at (-2,0) {Re$(L)$};
	
	%draw origins and label them
	\draw [fill] (0,0) circle [radius=0.05];
	\node [below] at (0,0) {$0$};
	\draw [fill] (-6,0) circle [radius=0.05];
	\node [above left] at (-6,0) {$0$};
	
	%Here I plot my unit circle
	%What I'm actually doing is plotting a bunch of 
	%quadrilaterals, each with its own color. I plot enough of them and I get the circle.
	\foreach \x in {0,0.0111,...,1} {
		\definecolor{currentcolor}{hsb}{\x, 1, 1}
		\draw[draw=none, fill=currentcolor]
		(360*\x - 2 : 0.96) -- (360*\x - 2 :1.04)
		-- (360*\x + 2:1.04) -- (360*\x + 2 :0.96) -- cycle;
	}
	
	%draw my final point, the one on the color wheel
	\draw [fill] (-0.8,-0.6) circle [radius=0.06];
	%label on the color wheel point
	%\node [below left] at (-0.8,-0.6) {$\mathcal{C}_{\mathbb{D},1}$};
	\end{tikzpicture}	
	\caption{Construction of the coloring map $\mathcal{C}_{\mathbb{D},1}$ for uniquely coloring \emph{asymptotic} parallel lines.}
	\label{fig:Poincare1}
\end{figure}

\begin{theorem}
	Let $f:\mathbb{C}_+\rightarrow \mathbb{C}_+$ represent a function on hyperbolic space. The function
	\begin{align*}
	\mathcal{C}_{\mathbb{D},1}:\mathbb{D}&\rightarrow [0,2\pi) \\
	z &\mapsto (\arg \circ T_{\mathbb{D}} \circ {\rm Re} \circ f \circ T_{\mathbb{D}}^{-1})(z)
	\end{align*} 
	assigns a unique color to the preimage of each tractrix line.
\end{theorem}
Figure~\ref{fig:Poincare1} both illustrates the map and serves as a visual proof. Notice that the result is still true if $T_{\mathbb{D}},T_{\mathbb{D}}^{-1}$ are replaced by simple inversion in the circle without the conjugation step, and we omit the conjugation step from Figure~\ref{fig:Poincare1} for readability.

\begin{figure}[h]
	\begin{subfigure}{0.35\linewidth}
		\begin{center}
			\vspace{1cm}
			\includegraphics[width=\linewidth]{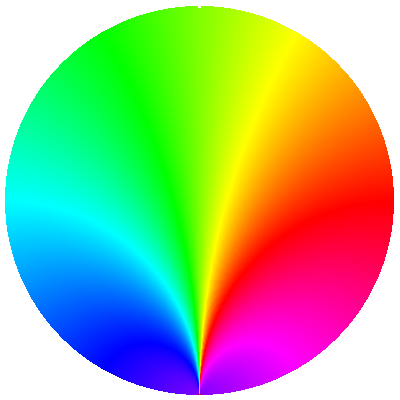}
		\end{center}
	\end{subfigure}
	\begin{subfigure}{0.63\linewidth}
		\begin{center}
			\includegraphics[width=\linewidth]{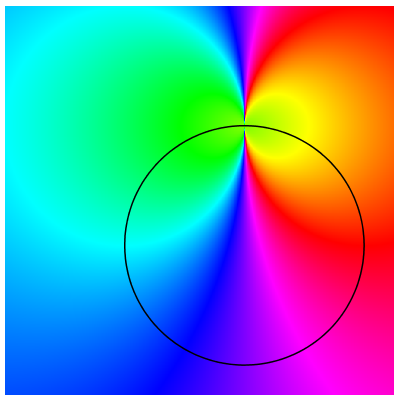}
		\end{center}
	\end{subfigure}
	\caption{The color assignment $\mathcal{C}_{\mathbb{D},1}$ for the Poincar\'{e} Disc is a \emph{true} phase portrait; at left is how it would appear with conjugation removed.}
	\label{fig:Poincarecolorings}
\end{figure}

The coloring is shown where $f$ is the identity at right in Figure~\ref{fig:Poincarecolorings}; at left is shown how it would appear without conjugation. Here it can be seen that this is a phase plot in the \emph{truest} sense of the word: the mapping returns to us the \emph{phase} of $T_{\mathbb{D}} \circ {\rm Re} \circ f \circ T_{\mathbb{D}}^{-1}$. 

Figure~\ref{fig:Poincarecolorings} also shows an important property of tractrix generators: they are \emph{asymptotic} parallel lines all meeting at the point at infinity corresponding to the top of the pseudosphere. Under the color map $\mathcal{C}_{\mathbb{D},1}$, we can always read off the preimage of this point at infinity as the point at infinity where all of the colors meet on the edge of the disc.

\begin{figure}[h]
	\begin{center}
		\vspace*{\fill}
		\includegraphics[width=.55\linewidth]{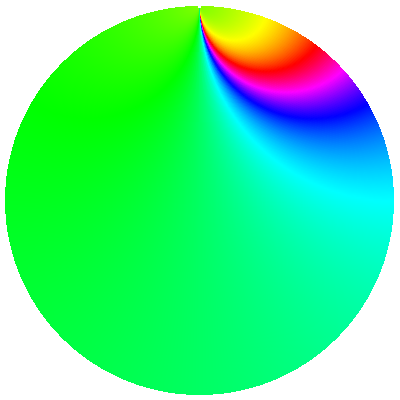}\vspace*{\fill} \includegraphics[width=.4\linewidth]{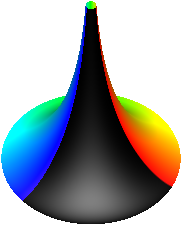}\vspace*{\fill}
	\end{center}
	\caption{A limit rotation on hyperbolic space.}
	\label{fig:poincarerotation}
\end{figure}

This unique property affords us insights about \emph{limit} rotations which we would not have obtained from our view on the pseudosphere alone. Letting $L_1,L_2$ be vertical lines with real parts $0,1$ respectively, we obtain a limit rotation which leaves fixed the point at infinity corresponding to the top of the pseudosphere. Its M\"{o}bius transformation $f$ is $z\mapsto z-2$. 

At left in Figure~\ref{fig:poincarerotation} we see the colors twist counterclockwise. Remembering that we are coloring points according to their images under $f$, this counterclockwise twisting really reveals a clockwise rotation of the space. However, in contradistinction to what we observed in the case of simple rotation, in the limit rotation we see that the point at infinity corresponding to the top of the pseudosphere remains fixed. This is because the tractrix lines from which we have constructed our direct motion intersect asymptotically at this point. 

At right in Figure~\ref{fig:poincarerotation}, the view on the pseudosphere reveals how tractrix lines are mapped to other tractrix lines. Additionally, it may be \emph{inferred} that the point at infinity corresponding to the top of the pseudosphere remains fixed, though it may not be seen directly.

\subsection{Poincar\'{e} Disc: Ultra-Parallel}

Part of what makes phase plotting so powerful in the complex case is its ability, combined with a modular plotting, to recover a full representation of the behavior of a function as in Figure~\ref{fig:complexmodular}. From a geometric standpoint, the height and color plot for a point $z$ indicate to the reader the unique half line and concentric circle about $0$ at whose intersection $f(z)$ lies.

To obtain a similar representative power in hyperbolic space requires finding a way to represent preimages for a different set of non-intersecting curves. Moreover, these curves must intersect each of the tractrix curves we chose for $\mathcal{C}_{\mathbb{D},1}$ no more than once.

The logical choice is another set of parallel lines. Since we already have a coloring map for asymptotic parallel lines, this time we opt for a set of \emph{ultra}-parallel lines, the semi-circular arcs centered at $0$ in $\mathbb{C}_+$.

Following the process we used to construct the color map $\mathcal{C}_{\mathbb{D},1}$, we seek a map which sends all of the points on a unique semicircle centered at $0$ to a unique real value. The modulus map $\mathcal{M}(z)=|z|$ will do. However, because this map sends unique ultra-parallel lines to only \emph{positive} real values, we cannot conclude with the map ${\rm arg}\circ T_\mathbb{D}$ as before. If we did so, we would end up with only values in $[0,\pi/2]\cup[3\pi/2,2\pi)$. Multiplication by $i$ before taking the argument serves to rotate the unit circle so that our image values return argument values in $[0,\pi]$. Squaring gives us back the full color wheel (equivalent to multiplying by 2 after applying the ${\rm arg}$ function).

\begin{figure}[h]
	\begin{tikzpicture}
	
	%draw the preimage points
	\draw [fill] (-5.13,-0.34) circle [radius=0.05];
	\draw [fill] (-4.63,-0.38) circle [radius=0.05];
	
	%draw the image points
	\draw [fill] (-5.56,1.92) circle [radius=0.05];
	\draw [fill] (-3.57,1.39) circle [radius=0.05];
	
	%draw T from disc to circle
	\draw [->,black, thick] (-5.1644,-0.1592) -- (-5.517,1.694);
	\node [below left] at (-5.517,1.694) {$T_{\mathbb{D}}^{-1}$};
	\draw [->,black, thick] (-4.5452,-0.2384) -- (-3.676,1.213);
	\node [left] at (-3.75,1.213) {$T_{\mathbb{D}}^{-1}$};
	
	%draw the points on the right semicircular arc
	\draw [fill] (-0.56,1.92) circle [radius=0.05];
	\draw [fill] (1.43,1.39) circle [radius=0.05];
	
	%draw the radius map point
	\draw [fill] (2,0) circle [radius=0.05];
	
	%arrows connecting right semicircular arc points with radius map point
	\draw [->, thick, gray] (-0.3808,1.7856) -- (1.744,.192);
	\draw [->,thick, gray] (1.5155,1.1815) -- (1.886,.278);
	\node [below left] at (1.6,1.1815) {$\mathcal{M}$};

	%draw the semicircular arc at right
	\draw [black,thick,domain=0:120] plot ({2*cos(\x)}, {2*sin(\x)});
	
	%draw the semicircular arc at left
	\draw [black,thick,domain=0:120] plot ({-5+2*cos(\x)}, {2*sin(\x)});

	%draw the squiggly function arrow
	\draw [->,thick,
	line join=round,
	decorate, decoration={
		zigzag,
		segment length=6,
		amplitude=2,post=lineto,
		post length=10pt				
	}] (-2.75,1.5) -- (-1.75,1.5);
	
	%label the squiggly function arrow
	\node [below] at (-2.25,1.4) {$f$}; 
	
	%draw the real axes
	\draw [gray] (-1.5,0) -- (2,0);
	\draw [gray] (-6,0) -- (-2.75,0);
	
	%Draw the left circle
	\draw [ultra thick] (-5,0) circle [radius=1];

	%draw origins and label them
	\draw [fill] (0,0) circle [radius=0.05];
	\node [below] at (0,0) {$0$};
	\draw [fill] (-5,0) circle [radius=0.05];
	\node [above] at (-5,0) {$0$};
	
	%Here I plot my unit circle
	%What I'm actually doing is plotting a bunch of 
	%quadrilaterals, each with its own color. I plot enough of them and I get the circle.
	\foreach \x in {0,0.0111,...,1} {
		\definecolor{currentcolor}{hsb}{\x, 1, 1}
		\draw[draw=none, fill=currentcolor]
		(360*\x - 2 : 0.96) -- (360*\x - 2 :1.04)
		-- (360*\x + 2:1.04) -- (360*\x + 2 :0.96) -- cycle;
	}
	
	%draw the image point on the color wheel
	\draw [fill] (.8,-0.6) circle [radius=0.06];
	
	%draw the final arrow to the point on the color wheel
	\draw [->,thick] (1.76,-0.12) -- (1.04,-0.48);
	\node [below right] at (1.4,-0.3) {$T_{\mathbb{D}}$};
	
	%point after multiplication by i
	\draw [fill] (.6,.8) circle [radius=0.05];
	
	%draw the arc for multiplication by i
	\draw [->,black,thick,domain=-25:45] plot ({.8*cos(\x)}, {.8*sin(\x)});
	\node [left] at (.625,.3) {$\cdot i$};

	\end{tikzpicture}	
	\caption{Construction of a phase plot $\mathcal{C}_{\mathbb{P},2}$ uniquely coloring \emph{ultra}-parallel lines on the Poincar\'{e} Disc.}
	\label{fig:Poincare2construction}
\end{figure}
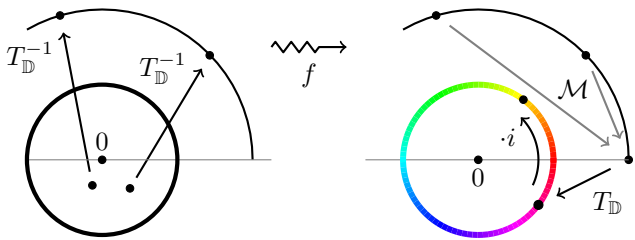

%\begin{figure}[ht]
%		\begin{center}
%		\includegraphics[width=.35\textwidth]{visualtheorem2}
%		\end{center}
%	\caption{Construction of a phase plot $\mathcal{C}_{\mathbb{P},2}$ uniquely coloring \emph{ultra}-parallel h-lines on the Poincar\'{e} Disc.}
%	\label{fig:Poincare2construction}
%\end{figure}

\begin{theorem}
	Let $f:\mathbb{C}_+\rightarrow \mathbb{C}_+$ represent a function on hyperbolic space. The function
	\begin{align*}
	\mathcal{C}_{\mathbb{D},2}:\mathbb{D}&\to [0,2\pi)\\
	z &\mapsto (2\cdot \arg \circ \;i\cdot T_{\mathbb{D}} \circ \mathcal{M} \circ f \circ T_{\mathbb{D}}^{-1})(z)
	\end{align*} 
	assigns a unique color to the preimage of each h-line corresponding to a circle in $\mathbb{C}_+$ centered at $0$.
\end{theorem}
Figure~\ref{fig:Poincare2construction} both illustrates the map and serves as a visual proof. Again the result remains true when $T_{\mathbb{D}},T_{\mathbb{D}}^{-1}$ are replaced by simple inversion in the circle without the conjugation step, and we again omit the conjugation step for the sake of clarity. This is \emph{also} a phase plot in the truest sense of the word, returning the \emph{phase} of $z \mapsto (i \circ T_{\mathbb{D}} \circ M \circ f \circ T_{\mathbb{D}}^{-1}(z))^2$, and so we again include the coloring in $\mathbb{C}$.

\begin{figure}[h]
	\begin{center} \includegraphics[width=\linewidth]{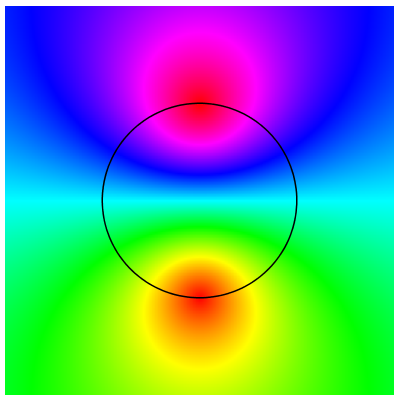}
	\end{center}
	\caption{The phase plot  $\mathcal{C}_{\mathbb{P},2}$.}
	\label{fig:Poincare2}
\end{figure}

The coloring where $f$ is the identity is shown in Figure~\ref{fig:Poincare2}. Interestingly, in this case the coloring without conjugation is exactly the same when $f$ is the identity. This difference between the two colorings is one consequence of the fact that $\mathcal{C}_{\mathbb{P},1}$ uniquely colors \emph{asymptotically} parallel lines while $\mathcal{C}_{\mathbb{P},2}$ uniquely colors \emph{ultra} parallel lines. To see why this difference is important, notice how colors along the edges of the discs in Figure~\ref{fig:Poincarecolorings} match everywhere except at the point at infinity on the rim where tractrix lines meet.

\begin{figure}[h]
	\begin{center}
		\vspace*{\fill}\includegraphics[width=.5\linewidth]{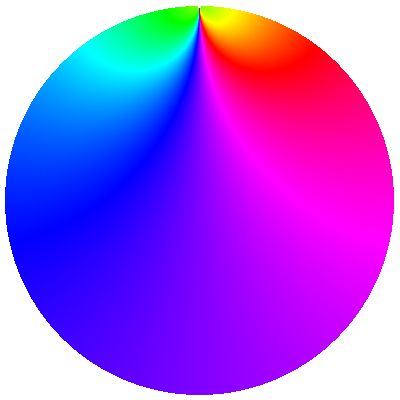}\vspace*{\fill}\includegraphics[width=.4\linewidth]{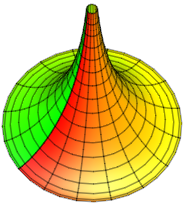}\vspace*{\fill}
	\end{center}
	\caption{A translation on hyperbolic space.}
	\label{fig:translationdown}
\end{figure}

With this new way of coloring plots on hyperbolic space, we consider the final direct motion: translation. Letting $L_1,L_2$ be semicircles centered at $0$ with radii $\sqrt3$ and $1$ respectively, we obtain the translation given by the M\"{o}bius transformation $f:z\mapsto z/3$. Plotting it in Figure~\ref{fig:translationdown} with $\mathcal{C}_{\mathbb{D},1}$ at left and $\mathcal{C}_{\mathbb{P}}$ at right, the behavior is not immediately clear.

\begin{figure}[h]
	\begin{center}
		\vspace*{\fill}\includegraphics[width=.55\linewidth]{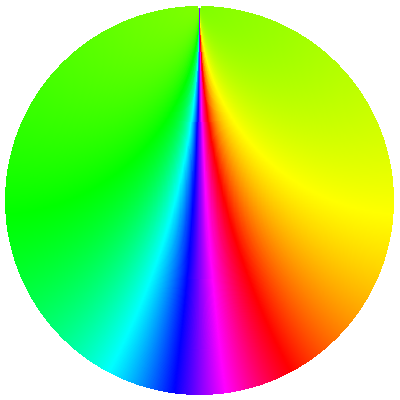}\vspace*{\fill}\includegraphics[width=.4\linewidth]{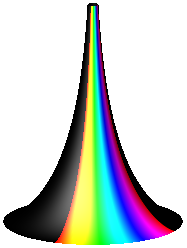}\vspace*{\fill}
	\end{center}
	\caption{A translation on hyperbolic space.}
	\label{fig:translationup}
\end{figure}

Letting $L_1,L_2$ be semicircles centered at $0$ with radii $1$ and $2$ respectively, we obtain the translation given by the M\"{o}bius transformation $f:z\mapsto 4z$. Plotting in Figure~\ref{fig:translationup} with $\mathcal{C}_{\mathbb{D},1}$ at left and $\mathcal{C}_{\mathbb{P}}$ at right, the behavior begins to come into focus. It appears that points may be sliding towards the top of the pseudosphere.

\begin{figure}[h]
	\begin{center}
		\vspace*{\fill}\includegraphics[width=.49\linewidth]{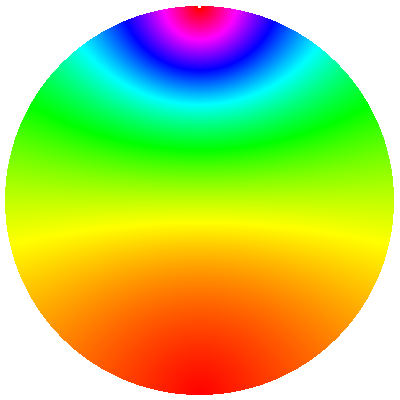}\vspace*{\fill}\includegraphics[width=.49\linewidth]{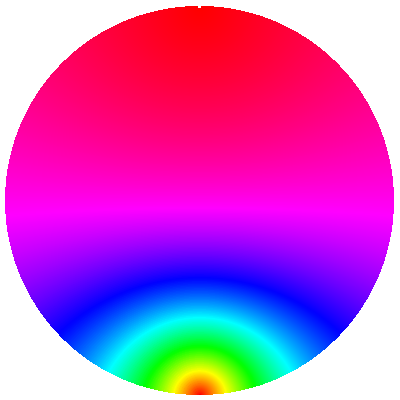}\vspace*{\fill}
	\end{center}
	\caption{Translations on hyperbolic space.}
	\label{fig:translation2ultra}
\end{figure}

Finally, using $\mathcal{C}_{\mathbb{D},2}$ to color $f:z\mapsto z/3$ at left and $f:z\mapsto 4z$ at right in Figure~\ref{fig:translation2ultra}, the behavior in both cases becomes clear, highlighting the advantage of having many views into hyperbolic space.

\begin{figure}[h]
	\begin{center}
		\vspace*{\fill}\includegraphics[width=.46\linewidth]{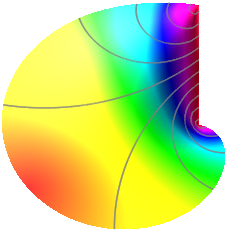}\vspace*{\fill}\includegraphics[width=.46\linewidth]{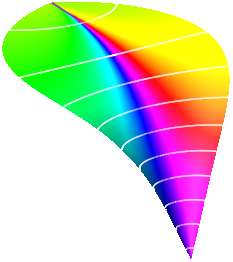}\vspace*{\fill}
	\end{center}
	\caption{Colored landscapes on $\mathbb{D}$.}
	\label{fig:poincaredynamic}
\end{figure}

With minor piecewise alterations we may combine both maps on $\mathbb{D}$ by using one for coloring and the other as a height map. At left in Figure~\ref{fig:poincaredynamic} we show the translation $f:z\mapsto z/3$ where $\mathcal{C}_{\mathbb{D},1}$ is the height map and $\mathcal{C}_{\mathbb{D},2}$ is the coloring map. At right we show the translation $f:z\mapsto 4z$ where $\mathcal{C}_{\mathbb{D},2}$ is the height map and $\mathcal{C}_{\mathbb{D},1}$ is the coloring map.

\begin{figure}[h]
	\begin{center}
		\vspace*{\fill}
		\adjustbox{trim={.21\width} {.21\height} {.21\width} {.21\height},clip}%
		{\includegraphics[width=.8\linewidth]{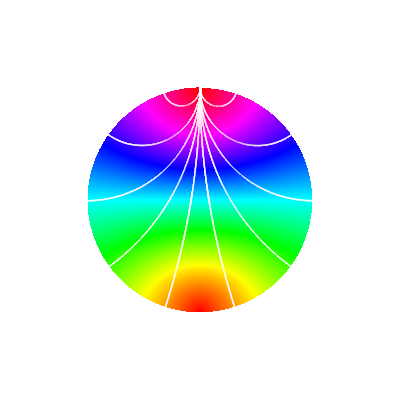}} \vspace*{\fill}
		\adjustbox{trim={.21\width} {.21\height} {.21\width} {.21\height},clip}%
		{\includegraphics[width=.8\linewidth]{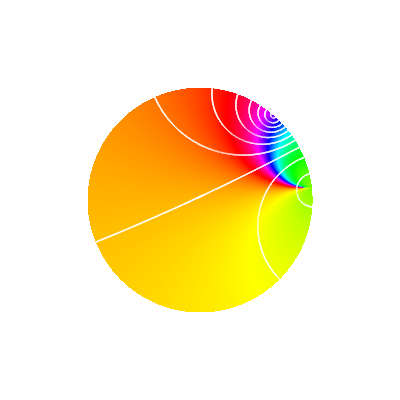}}\vspace*{\fill}
	\end{center}
	\caption{Colored contour plots on $\mathbb{D}$.}
	\label{fig:poincarecombined}
\end{figure}

The 3d plots in Figure~\ref{fig:poincaredynamic} may have the downside of distorting our view of $\mathbb{D}$, a problem easily remedied by using contours instead. We illustrate in Figure~\ref{fig:poincarecombined}; at left is the identity where $\mathcal{C}_{\mathbb{D},1}$ is the height map and $\mathcal{C}_{\mathbb{D},2}$ is the coloring map. At right we show the rotation corresponding to reflection in circles $L_1,L_2$ centered at $2$ and $-1$ with radii $13/8$ and $21/8$ respectively; $\mathcal{C}_{\mathbb{D},2}$ is the height map and $\mathcal{C}_{\mathbb{D},1}$ is the coloring map. The M\"{o}bius transformation is $f:z\mapsto (-249z-667)/(192z+215)$.

\subsection{Beltrami Half Sphere}

\begin{figure}[ht]
	\begin{subfigure}{0.6\linewidth}
	\begin{center}
		\includegraphics[width=\linewidth]{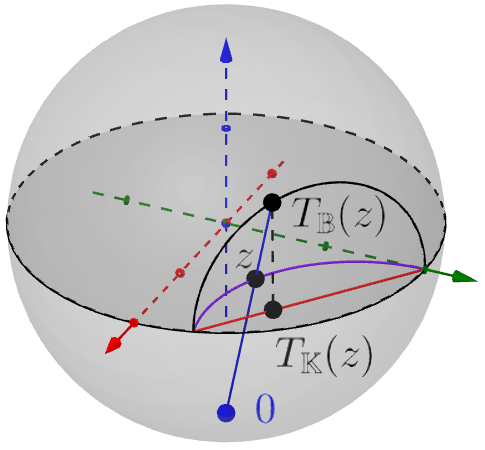}
	\end{center}
	\end{subfigure}
	\begin{subfigure}{0.38\linewidth}
		\begin{center}
		\adjustbox{trim={.21\width} {.21\height} {.21\width} {.3\height},clip}%
		{\includegraphics[width=1.6\linewidth]{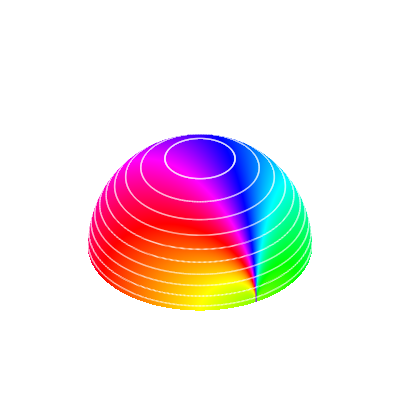}}\\
		\adjustbox{trim={.21\width} {.21\height} {.21\width} {.3\height},clip}%
		{\includegraphics[width=1.6\linewidth]{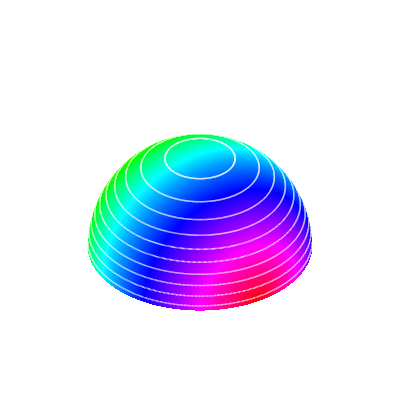}}
		\end{center}
	\end{subfigure}
	\caption{Construction of Beltrami half sphere and Klein Disc (left) and phase plotting on the Beltrami half sphere (right).}
	\label{fig:Beltrami}
\end{figure}

Having achieved our original goal of reinventing phase portraits to fully reveal the behavior of motions on hyperbolic space, we can consider the extension of our portraits to other representations. The \emph{Beltrami Half Sphere} $\mathbb{B}$ is obtained via \emph{lower} stereographic projection $T_{\mathbb{B}}$ from the Poincar\'{e} Disc $\mathbb{D}$ onto the upper half hemisphere with $\mathbb{D}$ as its base. This is illustrated at left in Figure~\ref{fig:Beltrami}. The natural coloring mappings for assigning unique colors to the preimages of unique lines are
\begin{align*}
\mathcal{C}_{\mathbb{B},1}:\mathbb{B}\rightarrow[0,2\pi)\quad\text{by}&\quad z\mapsto (\mathcal{C}_{\mathbb{D},1}\circ T_{\mathbb{B}}^{-1})(z)\\
\mathcal{C}_{\mathbb{B},2}:\mathbb{B}\rightarrow[0,2\pi)\quad\text{by}&\quad z \mapsto (\mathcal{C}_{\mathbb{D},2}\circ T_{\mathbb{B}}^{-1})(z).
\end{align*}

They are illustrated at right in Figure~\ref{fig:Beltrami} with  $\mathcal{C}_{\mathbb{B},1}$ at top and $\mathcal{C}_{\mathbb{B},2}$ at bottom.

\subsection{Klein Disc}

\begin{figure}[h]
	\begin{center}
		\vspace*{\fill}
		\includegraphics[width=.46\linewidth]{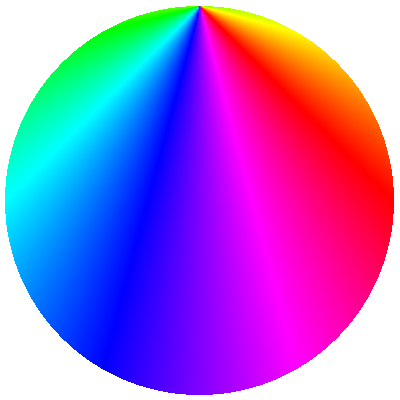}\vspace*{\fill}
		\includegraphics[width=.46\linewidth]{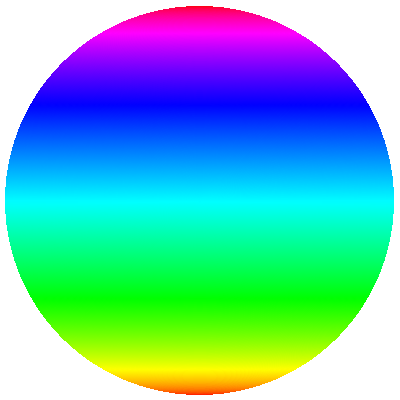}
		\vspace*{\fill}
	\end{center}
	\caption{Phase plotting on $\mathbb{K}$ with $\mathcal{C}_{\mathbb{K},1}$ (left) and $\mathcal{C}_{\mathbb{K},2}$ (right).}
	\label{fig:klein}
\end{figure}

The Klein disc $\mathbb{K}$ representation of hyperbolic space is obtained via the conformal map $T_{\mathbb{B}}$ to the Beltrami upper half sphere followed by vertical projection $T_{\mathbb{K}}$ which takes the upper half sphere to the disc. This is also shown in Figure~\ref{fig:Beltrami}. While this map is not conformal, this visualization has the interesting quality that h-lines are represented by lines we recognize in the Euclidean sense. Naturally the coloring mappings are
\begin{align*}
\mathcal{C}_{\mathbb{K},1}:\mathbb{K}\rightarrow[0,2\pi)\quad\text{by}&\quad z\mapsto (\mathcal{C}_{\mathbb{B},1}\circ T_{\mathbb{K}}^{-1})(z)\\
\mathcal{C}_{\mathbb{K},2}:\mathbb{K}\rightarrow[0,2\pi)\quad\text{by}&\quad z \mapsto (\mathcal{C}_{\mathbb{B},2}\circ T_{\mathbb{K}}^{-1})(z).
\end{align*}
These are shown in Figure~\ref{fig:klein}.

\subsection{Dini's Surface}

The limited view afforded by the pseudosphere cases like Figure~\ref{fig:translationdown} is an impediment we cannot completely work around. Still, we can improve this view.  

Motivated by our first example with a rotation, we use $\mathcal{C}_{\mathbb{P}}$ to plot a similar rotation on \emph{Dini's surface}, a twisted pseudosphere which affords us a view into hyperbolic space beyond the boundaries necessitated by the ordinary pseudosphere.

\begin{figure}[ht]
	\begin{center}
		\vspace*{\fill}
		\includegraphics[width=.49\linewidth]{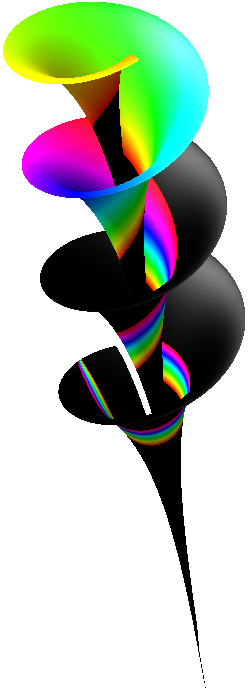}\vspace*{\fill}
		\includegraphics[width=.49\linewidth]{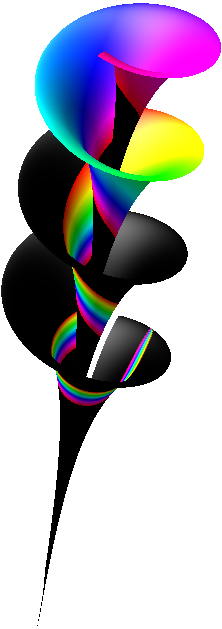}
		\vspace*{\fill}
	\end{center}
	\caption{Rotations on Dini's surface.}
	\label{fig:dinirotation}
\end{figure}
	
With an expanded view we naturally choose $L_1$ and $L_2$ to be larger circles, of radius $4\pi$, centered at $4\pi+\pi/8$ and $\pi/8$ respectively. The motion is plotted at right in Figure~\ref{fig:dinirotation} where we prune the surface below at $0$ and above at $7\pi$. The M\"{o}bius transformation is
$$f: z \mapsto (-261\pi z+1057\pi^2)/(-64z+8\pi).$$
As before, we see that the preimage of the point at infinity corresponding to the top of the pseudosphere has $\theta$ coordinate $\pi/8$, corresponding to the center of $L_2$. Geodesics curl around the extended surface, returning back to its edge near $\theta=4\pi+\pi/8$.

Here we can begin to have a bit of fun with symmetry. Let $L_1$ and $L_2$ have the same radius but instead be centered at $2\pi-(4\pi+\pi/8)$ and $2\pi-\pi/8$ respectively. Then the circles are reflections of their former selves about the tractrix line with real part $\pi$ which divides in half our set of tractrix lines whose preimages we color. The M\"{o}bius transformation is
$$f: z \mapsto (-136\pi z-769\pi^2)/(64z-120\pi).$$
If we twist the pseudosphere upwards from left to right---instead of downwards---and prune at $-5\pi$ above $2\pi$ below, we expect to see a mirror image of our first plot with the color spectrum reversed. In fact, this is what we do see at left in Figure~\ref{fig:dinirotation}.

\begin{figure}[ht]
	\begin{center}
		\includegraphics[width=.5\linewidth]{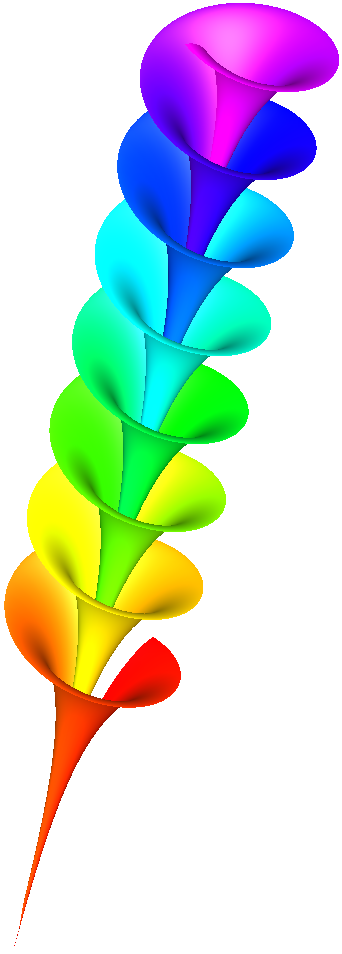}
	\end{center}
	\caption{Translation on Dini's surface.}
	\label{fig:dinitranslation}
\end{figure}

We conclude with a translation, a transformation which we think of as mundane in $\mathbb{C}$ but which comes alive in a thought-provoking way when phase plotted on $\mathbb{P}$ in Figure~\ref{fig:dinitranslation}.
	
Letting $L_1$ and $L_2$ be circles centered at $0$ with radii $3$ and $1$ respectively, we are forced to confront the fact that our intuitive notions of space and distance do not well prepare us for the infinite continuation of hyperbolic space beyond the pseudosphere's rim. We face this same challenge when considering the plot on $\mathbb{P}$ in Figure~\ref{fig:translationup}, although the two plots \emph{feel} quite different from one another because this translation is actually more similar to the one we had difficulty interpreting in Figure~\ref{fig:translationdown}. 

In \emph{this} case we are rescued by the augmented view afforded by Dini's surface. Pruning at $0$ and $15\pi$, we see that the entire visible region maps into the region with $\theta$ coordinates between $0$ and $2\pi$, and so the entire surface is painted with color. The M\"{o}bius transformation is $f: z\mapsto z/9$.

The use of visualization has already done much to change the way analysis is learned and taught. Phase portraits in particular have done much to illuminate the behavior of complex functions. We believe they can do the same for hyperbolic geometry. Moreover, as the popularity of the video \emph{M\"{o}bius Transformations Revealed} has demonstrated, mathematical beauty in full color appeals not only to mathematicians but also to a more general audience \cite{ARvideo}. Perhaps it is time we see more of it.\\

\textbf{Software}

The \emph{Maple} code used to generate all of the images in this paper is available at \url{https://ogma.newcastle.edu.au/vital/access/manager/Repository/uon:29853}.

%\textbf{Note}
%
%We extend our gratitude to Brailey Sims and Heinz Bauschke for their feedback and suggestions on our manuscripts. 

\textbf{Dedication}

This project is dedicated to the memory of Jonathan M. Borwein, our adviser, mentor, and friend.

%\section{Appendix: the Maple commands}

%Paul adds a brief appendix here which will provide a user guide to the code (will go online).


\begin{thebibliography}{99}

\bibitem{AR}D. N. Arnold and J. Rogness, ``M\"{o}bius transformations revealed.'' \emph{Notices of the AMS}, 55 (2008).

\bibitem{ARvideo}D. N. Arnold and J. Rogness, [jonathanrogness]. (2007, June 3). \emph{Moebius Transformations Revealed} [Video file]. Retrieved from \url{https://www.youtube.com/watch?v=JX3VmDgiFnY}

\bibitem{Tools}J.M. Borwein, ``The Life of Modern Homo Habilis Mathematicus: Experimental Computation and Visual Theorems.'' \emph{Tools and Mathematics}, 23-90, in \emph{Mathematics Education Library}, 347, Springer, 2016. 	
 	
 \bibitem{calendar}J.M. Borwein and A. Straub, ``Moment function of a 4-step planar random walk,'' \emph{Complex Beauties 2016}, (2016 Calendar). Available at: \url{http://www.mathe.tu-freiberg.de/files/information/calendar2016eng.pdf} 	
 	
\bibitem{Littlewood}J.E. Littlewood, \emph{A mathematician’s miscellany,} London: Methuen (1953); J. E. Littlewood and B\'{e}la Bollob\'{a}s, ed., \emph{Littlewood’s miscellany}, Camb. Univ. Press,
1986.
 	
\bibitem{Needham} T. Needham, \emph{Visual Complex Analysis}. Oxford University Press, 1997.

\bibitem{Wegert}E. Wegert, \emph{Visual Complex Functions: An Introduction with Phase Portraits}. Springer, 2012.

\bibitem{WS}E. Wegert and G. Semmler, ``Phase Plots of Complex Functions: A Journey in Illustration,'' \emph{Notices of the American Mathematical Society} \textbf{58(6)}, 2011.

\end{thebibliography}
\end{document}